%% file: complex_volumes6.tex
\documentclass[11pt]{article}

\usepackage{amsmath,amsfonts,amsthm,amssymb,graphicx, enumerate}
\usepackage[all,2cell,ps]{xy}

\theoremstyle{plain}
\newtheorem{thm}{Theorem}[section]
\newtheorem{lem}[thm]{Lemma}
\newtheorem{prop}[thm]{Proposition}
\newtheorem{cor}[thm]{Corollary}

\theoremstyle{definition}
\newtheorem{rmk}[thm]{Remark}

\newenvironment{pf}{\begin{proof}}{\end{proof}}

\DeclareMathOperator{\SL}{SL}
\DeclareMathOperator{\SU}{SU}
\DeclareMathOperator{\PU}{PU}
\DeclareMathOperator{\G}{G}
\DeclareMathOperator{\Res}{Res}
\newcommand{\ResI}{\Res^{(1)}}

\DeclareMathOperator{\tA}{A}
\DeclareMathOperator{\T}{T}

\DeclareMathOperator{\Disc}{D}
\DeclareMathOperator{\HH}{H}
\DeclareMathOperator{\CG}{C}
\DeclareMathOperator{\Aut}{Aut}
\DeclareMathOperator{\Ram}{Ram}
\DeclareMathOperator{\No}{N_{\ell / \Q}}
\newcommand{\Wo}{W^{\mathcal O}}
\DeclareMathOperator{\Id}{Id}

\newcommand{\Q}{\mathbb Q}
\newcommand{\R}{\mathbb R}
\newcommand{\Z}{\mathbb Z}
\newcommand{\C}{\mathbb C}
\newcommand{\A}{\mathbb A}
\newcommand{\F}{\mathbb F}

\newcommand{\Hysp}{\mathbf H}
\newcommand{\LL}{\mathbf L}

\newcommand{\q}{\mathfrak q}

\newcommand{\p}{\mathfrak p}
\newcommand{\TT}{\mathbf T}
\newcommand{\Gm}{\mathbf G_m}
\newcommand{\Op}{\mathcal O_{\mathfrak p}}

\newcommand{\cI}{\mathcal I}
\newcommand{\cP}{\mathcal P}
\newcommand{\cC}{\mathcal C}

\newcommand{\mvol}{\nu}
\newcommand{\EP}{\mu^{\mbox{\tiny EP}}}
\newcommand{\Lat}{\mathbf L}

\newcommand{\nb}{\mathfrak n}
\newcommand{\nbo}{\mathfrak n'}
\newcommand{\cln}{\mathfrak h}
\newcommand{\vol}{\mathrm{vol}}

\newcommand{\ind}{\epsilon_\ell}

\newcommand{\MM}{\overline{\mathrm{M}}}
\newcommand{\cMM}{\overline{\mathcal{M}}}

\newcommand{\Vf}{V_{\mathrm f}}
\newcommand{\ccol}{\mathbf P}
\newcommand{\CC}{\mathfrak C}
\newcommand{\bG}{\overline \G}
\newcommand{\Af}{\A_{\mathrm f}}
\newcommand{\bP}{\overline P}

\newcommand{\Idl}{I}
\newcommand{\Idlo}{\Idl^\emptyset}

\title{Covolumes of nonuniform lattices in $\PU(n, 1)$}

\author{Vincent Emery\footnote{Supported by SNSF projects  {\tt
PP00P2-128309/1} and {\tt PBFRP2-128067}}\\ \small{Universit\'e de
Gen\`eve}\\ \small{2-4 rue du Li\`evre, CP 64}\\ \small{1211 Gen\`eve 4, Switzerland}\\ \small{\textsf{vincent.emery@gmail.com}} \and Matthew Stover\footnote{Partially supported by NSF RTG grant DMS 0602191}\\ \small{University of Michigan}\\ \small{530 Church Street}\\ \small{Ann Arbor, MI 48109}\\ \small{\textsf{stoverm@umich.edu}}}

\date{\today}

\begin{document}

\maketitle

\begin{abstract} This paper studies the covolumes of nonuniform arithmetic lattices in $\PU(n, 1)$. We determine the smallest covolume nonuniform arithmetic lattices for each $n$, the number of minimal covolume lattices for each $n$, and study the growth of the minimal covolume as $n$ varies. In particular, there is a unique lattice (up to isomorphism) in $\PU(9, 1)$ of smallest Euler--Poincar\'e characteristic amongst all nonuniform arithmetic lattices in $\PU(n, 1)$. We also show that for each even $n$ there are arbitrarily large families of nonisomorphic maximal nonuniform lattices in $\PU(n, 1)$ of equal covolume.
\end{abstract}


\section{Introduction}\label{intro}

\input{intro6}


\section{Arithmetic subgroups of $\SU(n,1)$}\label{arithmetic}

\input{arithmetic6}


\section{Covolumes of nonuniform lattices}\label{nonuniform prasad}

\input{nonuniform6}


\section{Counting multiplicities}\label{multiplicities}

\input{multiplicities6}


\section{Theorems~\ref{thm:smallest overall}--\ref{thm:volume growth} and comparison with geometric volume estimates}\label{growth}

\input{growth6}


\bibliography{complex_volumes}

\end{document}

%% file: intro6.tex
In recent years, volumes of locally symmetric spaces have received a great deal of attention. One basic problem is to determine the smallest volume quotients of a given symmetric space of noncompact type by the action of a discrete group of isometries. A complete understanding of the smallest volume quotients of the hyperbolic plane follows from classical techniques. However, even for hyperbolic $3$-space a complete picture of the smallest volume orbifolds and manifolds was finished only very recently~\cite{Meyerhoff, Milley}. For progress in higher-dimensional hyperbolic spaces see~\cite{Hild}.

One observation is that all known smallest volume quotients of real hyperbolic $n$-space come from lattices that are arithmetically defined. The determination of the smallest volume arithmetic hyperbolic $n$-orbifolds is now complete~\cite{Chinburg--Friedman, Chinburg--Friedman--Jones--Reid, Belolipetsky, Emery, Belolipetsky--Emery}, and some progress has been made for complex hyperbolic space~\cite{Parker, Prasad--Yeung1, Prasad--Yeung2, Stover}. The purpose of this paper is to study the distribution of volumes of noncompact arithmetic quotients of complex hyperbolic space $\Hysp_\C^n$ and locate the smallest amongst them.

In other words, we study covolumes of nonuniform arithmetic lattices $\Gamma$ in $\PU(n, 1)$. Indeed, $\PU(n, 1)$ is the group of holomorphic isometries of $\Hysp_\C^n$ and the covolumes of $\Hysp_\C^n / \Gamma$ and $\PU(n, 1) / \Gamma$ are related by an explicit constant (see $\S$\ref{sec:overall}). Every nonuniform arithmetic subgroup $\Gamma$ of $\PU(n,1)$ is defined over $\Q$. More precisely, for each such $\Gamma$, there is an associated imaginary quadratic field $\ell = \Q\left(\sqrt{-d} \right)$ and a hermitian form on $\ell^{n + 1}$; see $\S$\ref{arithmetic} for the precise construction. We will say that $\Gamma$ is \emph{associated with $\ell / \Q$} to emphasize $\ell$, even though the algebraic group is defined rationally. 

We first introduce some notation necessary to state our results. For each imaginary quadratic field $\ell$, let $\Lat_\ell(n)$ be the set of isomorphism classes of nonuniform arithmetic lattices in $\PU(n, 1)$ associated with $\ell / \Q$. As normalization for the volume we use the Euler--Poincar\'e measure, i.e., the Haar measure $\EP$ on $\PU(n,1)$ such that $\EP(\PU(n,1)/\Gamma) = |\chi(\Gamma)|$ for any lattice $\Gamma$ in $\PU(n,1)$; here $\chi$ denotes the Euler--Poincar\'e characteristic in the sense of C.~T.~C.~Wall. Note that $\chi(\Gamma) = (-1)^n |\chi(\Gamma)|$ and the complex hyperbolic volume of $\Hysp_\C^n/\Gamma$ can be computed from $\chi(\Gamma)$ via the Chern--Gauss--Bonnet Theorem (see for instance~\cite{Hersonsky--Paulin}). Denote by $\mvol_\ell(n)$ the minimum of $|\chi(\Gamma)|$ among all lattices in $\Lat_\ell(n)$.  Furthermore, let $\nb_\ell(n)$ be the multiplicity of $\mvol_\ell(n)$, that is, the number of isomorphism classes in $\Lat_\ell(n)$ of covolume $\mvol_\ell(n)$. Note that since $\PU(n,1)$ has a nontrivial outer automorphism, a single element $\Gamma \in \Lat_\ell(n)$ determines two conjugacy classes of embeddings in $\PU(n,1)$ (corresponding to the two possible complex structures on the orbifold $\Hysp_\C^n/\Gamma$), but no more by Mostow--Prasad rigidity.

To give an exact formula for the covolume, we consider the $L$-function $L_\ell = \zeta_\ell / \zeta$ associated with $\ell/\Q$, where $\zeta_\ell$ is the Dedekind zeta function of $\ell$ and $\zeta = \zeta_\Q$ is the Riemann zeta function. Let $\Disc_\ell$ be the absolute value of the discriminant of $\ell$ and $r_\ell$ the number of finite primes ramified in $\ell / \Q$, i.e., the number of primes dividing $\Disc_\ell$. Finally, let $h_{\ell, n + 1}$ be the order of the subgroup of the class group of $\ell$ consisting of elements whose order divides $n + 1$.


\begin{thm}\label{thm:even}
Suppose that $n \ge 2$ is even. For each imaginary quadratic field $\ell$, all lattices in the set $\Lat_\ell(n)$ are commensurable. The minimal Euler--Poincar\'e characteristic among these lattices is given by
\begin{eqnarray*}
\mvol_\ell(n) &=& \frac{n + 1}{2^n h_{\ell, n + 1}} \;
\prod_{j = 1}^{n / 2} \zeta(1 - 2 j) L_\ell(-2 j).
\label{min-even}
\end{eqnarray*}
For the number $\nb_\ell(n)$ of lattices in $\Lat_\ell(n)$ realizing the minimal covolume $\mvol_\ell(n)$ we have $2^{r_\ell} \leq \nb_\ell(n) \leq 2^{r_\ell} \cln_\ell(n)$, where $\cln_\ell(n)$ is a divisor of $h_{\ell, n + 1}$.  In particular, if $h_{\ell,n+1} = 1$ we have $\nb_\ell(n) = 2^{r_\ell}$.
\end{thm}


It is a classical fact that for all $j > 0$, $\zeta(1 - 2 j)$ and $L_\ell(-2 j)$ are rational multiples of Bernoulli and generalized Bernoulli numbers. See \cite[$\S$16.6]{Ireland--Rosen}. Since $r_\ell$ is the number of primes dividing the discriminant of $\ell$, it becomes arbitrarily large as we vary $\ell$. Indeed, the discriminant of $\Q(\sqrt{-d})$ (with $d$ square free) equals $-d$ or $-4 d$ according to whether or not $d$ is congruent to $3$ modulo $4$. Therefore, we immediately obtain the following corollary.


\begin{cor}\label{cor:families-same-volume}
For all even $n$ there exist arbitrarily large families of nonisomorphic maximal lattices in $\PU(n,1)$ of the same covolume.
\end{cor}


Recall that Wang's Theorem~\cite{Wang} states that for each $n > 1$ there are only finitely many lattices in $\PU(n, 1)$ with covolume bounded above by any constant. Therefore, the covolumes in Corollary~\ref{cor:families-same-volume} necessarily grow with the size of the family. See~\cite{Emery2} for an analogue of this corollary for nonmaximal lattices in more general semisimple Lie groups. 

In odd dimensions, our results are slightly weaker than Theorem~\ref{thm:even}, due to the more delicate arithmetic of the lattices under consideration. In this setting we prove the following.


\begin{thm}\label{thm:odd}
Suppose that $n \ge 3$ is odd and that $\ell$ is imaginary quadratic. Then
\begin{eqnarray*}
\mvol_\ell(n) &=& \frac{(-1)^{\frac{n + 1}{2}} (n + 1) \ind(n)}{2^n h_{\ell, n + 1}}\;
\zeta(-n) \; \prod_{j = 1}^{(n - 1) / 2} \zeta(1 - 2 j) L_\ell(-2 j),
\label{nu-odd}
\end{eqnarray*}
where $\ind(n) \ge 2$ is a divisor of $2^{r_\ell}$. In particular $\ind(n) = 2$ if $\ell / \Q$ is ramified at exactly one rational prime.

Suppose that $\ell/\Q$ is ramified at exactly one prime. We have $\cln_\ell(n)/2 \leq \nb_\ell(n) \leq \cln_\ell(n)$, where $\cln_\ell(n)$ is a divisor of $h_{\ell, n + 1}$ if $8 \nmid n + 1$, and $\cln_\ell(n)$ is a divisor of $2 h_{\ell,n + 1}$ otherwise.
\end{thm}


One can again express this quantity in terms of Bernoulli and generalized Bernoulli numbers. Finally, using Theorems~\ref{thm:even} and~\ref{thm:odd} it is possible to completely determine the minimal covolume arithmetic lattices in any dimension. This gives the following theorem, which generalizes the case $n = 2$ studied in~\cite{Stover}.


\begin{thm}\label{thm:smallest volume} 
For all $n > 1$, the smallest covolume nonuniform arithmetic lattices in $\PU(n, 1)$ are associated with $\Q(\sqrt{-3})/\Q$. There are exactly two isomorphism classes of minimal covolume lattices when $n$ is even. There are at most two when $n \equiv 7 \pmod{8}$.  For all other $n$ there is exactly one.
\end{thm}


Note that for $\ell = \Q(\sqrt{-3})$ we have $h_\ell = r_\ell = 1$, so that in every dimension we get an explicit formula for the Euler--Poincar\'e  characteristic of the lattices with smallest covolume.

The case $n = 1$ is exceptional. Indeed, $\Q(i)$ determines the unique commensurability class of nonuniform arithmetic lattices in $\PU(1, 1)$, which is the conjugacy class of the modular group. There is exactly one orbifold of minimal volume, namely, the modular surface. We also note that only finitely many commensurability classes of nonarithmetic lattices in $\PU(n, 1)$ are known. Between the seven lattices constructed by Deligne and Mostow and the two new examples of Deraux--Parker--Paupert~\cite{Mostow, Deraux--Parker--Paupert}, there are between six and nine known commensurability classes of nonarithmetic lattices in dimension $2$. There is one known in dimension $3$~\cite{Mostow}, and there are no examples for higher $n$. These do not give lattices of smaller covolume than the minimal covolume arithmetic lattices (see the remark in \cite{Stover}).

We now briefly describe the strategy for proving the above results. The first step is to consider Prasad's volume formula \cite{Prasad}, which gives the covolume of a so-called principal arithmetic lattice. It is known that any maximal lattice is the normalizer of a principal arithmetic lattice, and the index is determined by analysis of Galois cohomology, particularly the results of Borel and Prasad \cite{Borel--Prasad}. Our results on the number of lattices realizing the minimal covolume requires further close analysis of automorphisms of algebraic groups and Galois cohomology.

We close the paper in $\S$\ref{growth} by studying the behavior of the smallest covolume as $n$ varies. We prove that the smallest Euler--Poincar\'e characteristic (in absolute value) is realized by the unique lattice of smallest covolume in $\Lat_{\Q(\sqrt{-3})}(9)$. Let $\Gamma_0$ denote this lattice. It is commensurable with the lattice acting on $\mathbf H_\C^9$ constructed by Deligne--Mostow \cite{Mostow}. This Deligne--Mostow lattice also appears in work of Daniel Allcock \cite{Allcock}. It turns out that $\Gamma_0$ also has smallest overall covolume with respect to the complex hyperbolic metric of constant biholomorphic curvature $-1$.


\begin{thm}\label{thm:smallest overall}
Let $\Gamma_0 < \PU(9,1)$ be as above, and normalize the holomorphic sectional curvatures of $\Hysp_\C^n$ to $-1$. For any $n \geq 2$, let $\Gamma < \PU(n, 1)$ be a nonuniform arithmetic lattice. Then $|\chi(\Gamma)| \geq |\chi(\Gamma_0)|$ and $\vol(\mathbf H_\C^n / \Gamma) \geq \vol(\mathbf H_\C^9 / \Gamma_0)$ with equality if and only if $\Gamma$ is isomorphic to $\Gamma_0$. The Euler--Poincar\'e characteristic and complex hyperbolic covolume of $\Gamma_0$ are:
\begin{eqnarray*}
\chi(\Gamma_0) &=& \frac{-809}{5,746,705,367,040} \approx -1.40776
\times 10^{-10} \, ;\\
\vol(\mathbf H_\C^9 / \Gamma_0) &=&  \frac{809 \pi^9}{79,550,340,408,000}
\approx 3.03148 \times 10^{-7} \, .
\end{eqnarray*}
\end{thm}


We also determine the asymptotic behavior of the minimal covolume and compare it to known volume bounds for complex hyperbolic $n$-orbifolds. In particular, we prove the following.


\begin{thm}\label{thm:volume growth}
As $n \to \infty$, the ratio $\q_\ell(n) = \mvol_\ell(n+1) / \mvol_\ell(n)$ grows super-exponentially. That is, for each $\alpha > 0$ there exist $n_0 \in \mathbb{N}$ such that $\q_\ell(n) > e^{\alpha n}$ for all $n \ge n_0$. In particular, $\mvol_\ell(n)$ grows super-exponentially.
\end{thm}


\subsubsection*{Acknowledgments}

This work was started at the Oberwolfach Seminar `On arithmetically defined hyperbolic manifolds' and was completed during the London Mathematical Society EPSRC Durham Symposium `Geometry and arithmetic of lattices'. We thank the organizers of both of these excellent workshops. We also thank Gopal Prasad for conversations related to this work and the referees for many helpful suggestions.


%% file: arithmetic6.tex
While our goal is to analyze covolumes of lattices in $\PU(n, 1)$, many results on the structure of arithmetic subgroups require one to consider lattices in a simply connected algebraic group. Therefore, we must spend much of this paper instead considering lattices in $\SU(n, 1)$. We later return to $\PU(n, 1)$ via the $(n + 1)$-fold covering $\SU(n, 1) \to \PU(n, 1)$.


\subsection{Hermitian forms and lattices in $\SU(n, 1)$}\label{hermitian forms and lattices}


Let $\G$ be a semisimple $\Q$-algebraic group. Then $\G(\Z)$ is a lattice in the real Lie group $\G(\R)$~\cite{Borel--Harish-Chandra}. Any $\Gamma < \G(\R)$ commensurable with $\G(\Z)$ is called an \emph{arithmetic lattice} in $\G(\R)$. The algebraic groups of interest for us are isotropic of type ${}^2 \tA_n$~\cite{Tits}. All such groups are constructed as follows (cf., \cite[$\S$4.4]{Prasad--Yeung2}).


\begin{prop}\label{defined over Q}
Let $\Gamma < \SU(n, 1)$ be a nonuniform arithmetic lattice. Then, there exists an imaginary quadratic field $\ell$ and a hermitian form $h$ of signature $(n, 1)$ on $\ell^{n + 1}$ such that $\Gamma$ is commensurable with $\G(\Z)$, where $\G$ is the $\Q$-algebraic group $\SU(h)$ of special unitary automorphisms of $h$.
\end{prop}


We now describe how the classification of unitary groups over nonarchimedean local fields determines the commensurability class of an arithmetic lattice. If $p$ is a prime that splits in $\ell$, i.e., a prime that splits as a product of two distinct prime ideals in $\ell$, then $\SU(h)$ is isomorphic over $\Q_p$ to $\SL_{n + 1}(\Q_p)$. If $p$ does not split, let $\mathfrak p$ be the unique prime ideal of $\ell$ over $p$. Then the completion $\ell_{\mathfrak p}$ of $\ell$ at $\mathfrak p$ is a quadratic extension of $\Q_p$ and $\G(\Q_p)$ is the unitary group of $h$ considered as a hermitian form on $\ell_{\mathfrak p}^{n + 1}$. For a nonsplit place, the following classifies the possible nonarchimedean completions of $\G$.


\begin{prop}\label{classification of p-adic groups}
Let $L / \Q_p$ be a quadratic extension, $h_1, h_2$ be hermitian forms on $L^{n + 1}$, and $\G_1, \G_2$ the associated $\Q_p$-algebraic groups. If $n$ is even, then $\G_1(\Q_p) \cong \G_2(\Q_p)$ and they are quasi-split. If $n$ is odd, then $\G_1(\Q_p) \cong \G_2(\Q_p)$ if and only if $\det(h_1)$ is congruent to $\det(h_2)$ modulo $\mathrm N_{L / \Q_p}(L^\times)$ and $\G_j$ is quasi-split if and only if $\det(h_j)$ is a norm from $L^\times$.
\end{prop}


It follows from the isomorphism between $\mathrm N_{L / \Q_p}(L^\times) / \Q_p^\times$ and $\mathrm{Gal}(L / \Q_p)$ in elementary local class field theory~\cite[Ch.~5]{Swinnerton-Dyer} that there is exactly one isomorphism class when $n$ is even and there are exactly two when $n$ is odd. This has the following consequence for the classification of nonuniform arithmetic lattices in $\SU(n, 1)$.


\begin{cor}\label{classification of nonuniform lattices}
If $n$ is even, associated with each field $\ell$ there is exactly one commensurability class of nonuniform arithmetic lattices in $\SU(n, 1)$. When $n$ is odd, commensurability classes of nonuniform arithmetic lattices in $\SU(n, 1)$ are in one-to-one correspondence with elements of $\Q^\times / \No(\ell^\times)$.
\end{cor}


Let $\ell$ be an imaginary quadratic field, $h$ a hermitian form on $\ell^{n + 1}$, and $\T$ be the finite set of (necessarily non-split) places $p$ of $\Q$ where $\det(h)$ fails to be a local norm from $\ell_{\p}^\times$, where $\mathfrak p$ is the unique prime of $\ell$ above $p$. Let $\G_{\ell, \T}$ be the associated $\Q$-algebraic group. Then $\G_{\ell, \T}$ fails to be quasi-split precisely at the places in $\T$.

When $\T = \emptyset$, which is always the case for $n$ even, we write $\G_\ell$ instead. Then $\G_\ell$ is the unique isotropic $\Q$-algebraic group with the following properties:
\begin{enumerate}[(i)]

\item $\G_\ell(\R) \cong \SU(n, 1)$;

\item $\G_\ell$ becomes an inner form over $\ell$;

\item $\G_\ell$ is quasi-split over $\Q_p$ for every rational prime $p$.

\end{enumerate}
For $\T \neq \emptyset$, $\G_{\ell, \T}$ satisfies (i) and (ii), but not (iii).

Now, we refine our description of the local structure of $\G_{\ell, \T}$. This will be crucial for parametrizing and studying the covolumes of the various maximal lattices in a commensurability class. This local structure can be understood by the \emph{local index} attached to semisimple groups over local fields, following Bruhat--Tits theory. We will use the notation and classification given in~\cite[$\S$4]{Tits}.

We begin with $\G_\ell$, i.e., the case $\T = \emptyset$. If $p$ is unramified in $\ell / \Q$ the local  index of $\G_\ell / \Q_p$ is $\tA_n$ or ${}^2 \tA^\prime_n$ according to whether or not $\G_\ell / \Q_p$ splits.  Note that the local index ${}^2 \tA^\prime_n$ differs with the parity of $n$. Now, suppose that $p$ is ramified in $\ell / \Q$. Then $\G_\ell / \Q_p$ has local index $\mbox{C--BC}_m$ if $n = 2 m$ and $\mbox{B--C}_m$ if $n = 2 m - 1$ (resp.~$\mbox{C--B}_2$ for $n = 3$). The latter index has a nontrivial automorphism, whereas the former is symmetry-free.
 
Now consider the case $\T \neq \emptyset$. When $\G_{\ell, \T} / \Q_p$ is quasi-split, its local structure is similar to $\G_{\ell, \T} / \Q_p$. If $\G_{\ell, \T} / \Q_p$ is not quasi-split it has local index ${}^2 \mbox{B--C}_m$ or ${}^2 \tA^{\prime \prime}_{2 m - 1}$  according to whether or not $p$ ramifies in $\ell / \Q$.


\subsection{Principal arithmetic lattices}\label{principal arithmetic}


Let $\ell$ be an imaginary quadratic field and $\G_{\ell, \T}$ the $\Q$-algebraic group associated with a hermitian form on $\ell^{n + 1}$ as above. We now describe the principal arithmetic lattices in $\G_{\ell, \T}(\Q)$.

Let $\Vf$ be the set of finite primes of $\Q$. For each $p \in \Vf$, let $P_p$ be a parahoric subgroup of $\G_{\ell, \T}(\Q_p)$~\cite[$\S$3.1]{Tits2}, and set $P = (P_p)_{p \in \Vf}$. We assume that $P$ is \emph{coherent}, so the restricted direct product $\prod_p P_p$ defines an open compact subgroup of $\G_{\ell, \T}(\Af)$, where $\Af$ are the finite adeles of $\Q$.

Embed $\G_{\ell, \T}(\Q)$ in $\G_{\ell, \T}(\Af)$ diagonally. By construction the subgroup
\[
\Lambda_{\ell, \T}(P) = \G_{\ell, \T}(\Q) \cap \prod_{p \in \Vf} P_p
\]
of $\G_{\ell, \T}(\R) \cong \SU(n, 1)$ is commensurable with $\G_{\ell, \T}(\Z)$. These are the so-called \emph{principal} arithmetic lattices.

%% file: nonuniform6.tex

\subsection{Covolumes of principal arithmetic lattices}\label{sec:constuction-of-Lambda}


We retain all notation from $\S$\ref{arithmetic}. In particular, $\ell$ is an imaginary quadratic field, $\G_{\ell, \T}$ is the $\Q$-algebraic group determined by a hermitian form on $\ell^{n + 1}$ such that $\G_{\ell, \T}$ fails to be quasi-split precisely at the rational primes in $\T$, and $\G_\ell$ denotes the unique such group that is quasi-split at every rational prime $p$ (i.e., where $\T = \emptyset$).

Let $\ccol_{\ell,\T}$ (resp., $\ccol_\ell$ if $\T = \emptyset$) be the set of coherent collections $P = (P_p)_{p \in \Vf}$ such that each parahoric subgroup $P_p$ is of maximal volume in $\G_{\ell,\T}(\Q_p)$. Then, for $P \in \ccol_{\ell,\T}$, all subgroups $P_p$ are special~\cite[A.5]{Borel--Prasad}. Moreover, if $p \not \in \T$  is not ramified in $\ell/\Q$, then $P_p$ is necessarily hyperspecial~\cite[\S 3.8]{Tits2}.
For $P \in \ccol_{\ell,\T}$ we consider the principal arithmetic
subgroup $\Lambda_{\ell,\T} = \Lambda_{\ell,\T}(P)$
(see $\S$\ref{principal arithmetic}). When $\T = \emptyset$, set $\Lambda_\ell = \Lambda_{\ell,\T}$.

The subgroup $\Lambda_{\ell,\T}$ depends on the choice of $P \in \ccol_{\ell,\T}$. However, by Prasad's volume formula, its covolume does not depend on $P \in \ccol_{\ell,\T}$. Let $s$ be $n(n + 3) / 4$ when $n$ is even and $(n - 1)(n + 2) / 4$ when $n$ is odd. For $\mu$ normalized as in~\cite{Prasad} we get that $\mu( \SU(n, 1) / \Lambda_\ell)$ equals:
\begin{eqnarray}
\label{eq:covol-Lambda}
\hspace{0.7cm}
\Disc_\ell^s \prod_{j = 1}^n \frac{j!}{(2 \pi)^{j + 1}} \;
\zeta(2) L_\ell(3) \cdots \left\{ \begin{array}{lr}
\hspace{-0.2cm}	\zeta(n + 1) & \mbox{$n$ odd}\\
\hspace{-0.2cm}	L_\ell(n + 1)&  \mbox{$n$ even}
\end{array} \right.
\label{vol-gL_p}
\end{eqnarray}
where $L_\ell$ is the Dirichlet $L$-function attached to the quadratic extension $\ell / \Q$ and $\Disc_\ell$ is the absolute value of the discriminant of $\ell$. We relate this to the Euler--Poincar\'e measure described in our introduction in $\S$\ref{sec:G-prime} below. See also~\cite[$\S$4]{Borel--Prasad}.

The covolume of $\Lambda_{\ell, \T}$ differs from the covolume of $\Lambda_\ell$ by a product of rational factors $\lambda_p > 1$ ($p \in \T$):
\begin{eqnarray}\label{volume-Lambda-prime}
\mu \left( \SU(n, 1) / \Lambda_{\ell, \T} \right)  &=& \mu \left( \SU(n, 1) / \Lambda_\ell \right)\; \prod_{p \in \T} \lambda_p.
\end{eqnarray}
More precisely, each factor $\lambda_p$ is given by 
\begin{eqnarray}
	\lambda_p &=& \frac{p^{\left( \dim \MM_p + \dim\cMM_p \right)/2}}{|\MM_p(\F_p)|}\; \frac{|\cMM_p(\F_p)|}{p^{\dim \cMM_p}},
	\label{eq:lambda-factor}
\end{eqnarray}
where the reductive $\F_p$-groups $\MM_p$ and $\cMM_p$ are defined as follows. Let $\overline{\G}_p$ (resp., $\overline{\mathcal G}_p$) be the $\F_p$-group associated with $\G_{\ell, \T}$ (resp., the group of which $\G_{\ell, \T}$ is an inner form) as in~\cite[$\S$3.5]{Tits2} and $R(\overline{\G}_p)$ (resp., $R(\overline{\mathcal G}_p)$) be the unipotent radical. Then $\MM_p$ (resp., $\cMM_p$) is the reductive $\F_p$-subgroup so that
\[
\overline{\G}_p = \MM_p R(\overline{\G}_p) \quad \left( \textrm{resp.,}\ \overline{\mathcal G}_p = \cMM_p R(\overline{\mathcal G}_p) \right).
\]
See~\cite[\S 2.2]{Prasad} for more details.

Note that $\Lambda_{\ell, \T}$ (including $\T = \emptyset$) is of minimal covolume among \emph{principal} arithmetic subgroups of $\G_{\ell, \T}$. However, as we shall see, $\Lambda_{\ell, \T}$ is usually not maximal in $\SU(n, 1)$. The normalizer $\Gamma_{\ell, \T} = N_{\SU(n, 1)}(\Lambda_{\ell, \T})$ is a maximal arithmetic subgroup of $\SU(n,1)$ (as usual, we write $\Gamma_\ell$ when $\T =\emptyset$). In fact, results from~\cite{Borel--Prasad} can be used to show that the indices $[\Gamma_{\ell, \T} : \Lambda_{\ell, \T}]$ do not depend on the choice of the coherent collections $P \in \ccol_{\ell,\T}$ and that $\Gamma_{\ell, \T}$ realizes the minimal covolume in its commensurability class (see~\cite[\S 12.3]{Emery} for details).


\subsection{Index computations for normalizers of principal arithmetic lattices}\label{index-comput}


See~\cite[Ch.\ 12]{Emery} for an explanation of how to estimate the indices $[\Gamma_{\ell, \T} : \Lambda_{\ell, \T}]$ based on the work of Borel--Prasad~\cite{Borel--Prasad}. The computation of these indices involves the center of the group $\G_{\ell, \T}$, which we denote by $\CG_\ell$. Let $\mu_{n + 1}$ be the $\ell$-group of units of order $n + 1$. The center of $\G_{\ell, \T}$ does not depend on $\T$ and is $\Q$-isomorphic to the group $\ResI_{\ell / \Q}(\mu_{n + 1})$ of norm $1$ elements in the restriction of scalars from $\ell$ to $\Q$.

We first restrict our attention to the computation of $[\Gamma_\ell:\Lambda_\ell]$. Until the end of $\S$\ref{ss:covolume-of-Gamma_ell} we consider the algebraic group $\G_\ell$. In particular, $\G_\ell/\Q_p$ is always quasi-split.

For each finite prime $p$, the Galois cohomology group $\HH^1(\Q_p, \CG_\ell)$ acts on the local Dynkin diagram $\Delta_p$ of $\G_\ell(\Q_p)$. The corresponding homomorphism is denoted $\xi_p : \HH^1(\Q_p, \CG_\ell) \to \Aut(\Delta_p)$. The map $\HH^1(\Q, \CG_\ell) \to \HH^1(\Q_p, \CG_\ell)$ determines a natural diagonal homomorphism
\begin{eqnarray}\label{xi}
\xi &:& \HH^1(\Q, \CG_\ell) \to \prod_p \Aut(\Delta_p).
\end{eqnarray}
Let $A = \ker \left( \HH^1(\Q, \CG_\ell)\to \HH^1(\R, \CG_\ell) \right)$ and $A_\xi$ be the kernel of $\xi$ restricted to $A$. We then have the following.


\begin{lem}\label{lem:index-Rholfs}
Let $\Lambda_{\ell}$ and $\Gamma_{\ell}$ be as above. Then
\begin{eqnarray}\label{index-Rholfs}
[\Gamma_\ell : \Lambda_\ell] &=& \frac{n + 1}{\left| \mu_{n + 1}(\ell)\right|} \;
\cdot\;
\left| A_\xi \right|.
\end{eqnarray}
\end{lem}



\begin{pf}
By~\cite[Prop. 2.9]{Borel--Prasad} and~\cite[Lemme 12.1]{Emery} we have
\begin{eqnarray}
[\Gamma_\ell : \Lambda_\ell] &=& \left| \CG_\ell(\R) / \CG_\ell(\Q) \right| \;
\cdot\;
\left| A_\xi \right|.
\end{eqnarray}

The group $\CG_\ell(\Q) \cong \ResI_{\ell / \Q}(\mu_{n + 1})(\Q)$ corresponds to the set of elements $x \in \mu_{n + 1}(\ell)$ for which $\No(x) = 1$. Since $\ell / \Q$ is imaginary quadratic, we have $\No \left( \mu_{n + 1}(\ell) \right) = 1$ and hence $\CG_\ell(\Q) \cong \mu_{n + 1}(\ell)$. The same argument shows that $\CG_\ell(\R)$ is isomorphic to $\mu_{n + 1}(\C)$ and hence has order $n + 1$. Then
\begin{eqnarray}\label{index-centre}
\left| \CG_\ell(\R) / \CG_\ell(\Q) \right| &=& \frac{n + 1}{\left| \mu_{n + 1}(\ell)\right|}.
\end{eqnarray}
This proves the lemma.
\end{pf}



\subsection{Actions on $\Delta_p$}\label{sec:A_xi_as_sbgp_of_LL}


Define the following subgroup of $\ell^\times$:
\begin{eqnarray}\label{LL}
\LL &=&  \left\{x \in \ell^\times \; / \; N_{\ell / \Q}(x) \in (\Q^\times)^{n + 1} \right\}.
\end{eqnarray}
The first Galois cohomology groups $\HH^1(\Q, \CG_\ell)$ and $\HH^1(\R, \CG_\ell)$ are described by following diagram, with exact rows:
\[
\xymatrix{
1 \ar[r] & \mu_{n + 1}(\Q) \ar[r] \ar[d]^{\cong} & \HH^1(\Q, \CG_\ell) \ar[d] \ar[r] & \LL / (\ell^\times)^{n + 1} \ar[r] \ar[d] & 1 \\
1 \ar[r] & \mu_{n + 1}(\R) \ar[r] &  \HH^1(\R, \CG_\ell) \ar[r] & 1 \ar[r] & 1
}
\]
From this we see that the first row splits and the subgroup $A \subset \HH^1(\Q, \CG_\ell)$ appearing in $\S$\ref{index-comput} can be identified with $\LL/(\ell^\times)^{n + 1}$. Then we get an action of $\LL$ on every diagram $\Delta_p$. The kernel $A_\xi$ fits into the following exact sequence:
\begin{eqnarray}\label{A_xi:exact-seq}
1 \quad \to \quad A_\xi \quad \to \quad \LL/(\ell^\times)^{n+1} \quad \to \quad \prod_p \Aut(\Delta_p) \quad \to \quad 1.
\end{eqnarray}

Now suppose that $p$ is a finite prime. The group $\HH^1(\Q_p, \CG_\ell)$ is described by the exact
sequence:
\begin{eqnarray}\label{exact-seq-H1Qp}
1 \quad \to \quad N_p \quad \to \quad \HH^1(\Q_p, \CG_\ell) \quad \to \quad W_p \quad \to 1,
\end{eqnarray}
where 
\begin{eqnarray}
N_p &=& \mu_{n + 1}(\Q_p) / \No(\mu_{n + 1}(\ell \otimes \Q_p)), \label{N_p} \\
W_p &=& \ker \left( (\ell \otimes \Q_p)^\times / ( (\ell \otimes \Q_p)^\times)^{n + 1} \stackrel{\No}{\longrightarrow} \Q_p^\times/(\Q_p^\times)^{n+1} \right). \label{W_p}
\end{eqnarray}
Recall that $\ell \otimes \Q_p$ is isomorphic to $\prod_{\p|p}\ell_\p$, which is a field if $p$ is either inert or ramified in $\ell/\Q$, or equal to two copies of $\Q_p$ if $p = \p \overline{\p}$ splits. In the latter case the norm $\No$ on $\ell \otimes \Q_p \cong \ell_\p \times \ell_{\overline{\p}}$ is obtained by multiplying the two components together, and the group $W_p$ is isomorphic to $\Q_p^\times/(\Q_p^\times)^{n + 1}$.


\begin{lem}\label{lem:N_p-trivial}
The subgroup $N_p \subset \HH^1(\Q_p,\CG_\ell)$ acts trivially on $\Delta_p$ for each prime $p$, i.e., $N_p \subset \HH^1(\Q_p, \CG_\ell)_{\xi_p}$.
\end{lem}



\begin{pf}
A proof for $p$ not ramified in $\ell / \Q$ is given in~\cite[\S 5.3]{Borel--Prasad}. Suppose that $p$ is ramified in $\ell / \Q$. If $n = 2 m$ is even the local index of $\G_\ell / \Q_p$ is $\mbox{C--BC}_m$, for which we have $\Aut(\Delta_p) = 1$. So in this case the lemma trivially holds.

Assume that $n$ is odd. Let $\TT$ be the $\Q_p$-group defined as the kernel of the map
\[
\begin{array}{ccc}
\Gm \times \Res_{\ell / \Q}(\Gm) &\to& \Gm \\
(x,y) &\mapsto& x^{n + 1} \No(y)^{-1},
\end{array}
\] 
where $\Gm$ denotes the multiplicative group. We get a commutative diagram with exact rows:
\[
\xymatrix{
1 \ar[r] & \CG_\ell \ar[r] \ar@{=}[d] & \Res_{\ell / \Q}(\mu_{n + 1}) \ar@{^{(}->}[d] \ar[r]^-{\No} & \mu_{n + 1} \ar[r] \ar@{^{(}->}[d] & 1 \\
1 \ar[r] & \CG_\ell \ar[r] &  \Res_{\ell / \Q}(\Gm) \ar[r]^-{\theta} & \TT \ar[r] & 1,
}
\]
the homomorphism $\theta$ being $\theta(z) = (\No(z), z^{n + 1})$, and the embedding $\mu_{n + 1} \hookrightarrow \TT$ by $x \mapsto (x, 1)$. Passing to cohomology, the first row gives the sequence~\eqref{exact-seq-H1Qp}, whereas the second gives an isomorphism between $\HH^1(\Q_p, \CG_\ell)$ and a quotient of
\begin{eqnarray*}
\TT(\Q_p) &=& \left\{ (x, y) \in \Q_p^\times \times \ell_\p^\times \;|\; x^{n + 1} = \No(y) \right\}.
\end{eqnarray*}
We refer to~\cite[Prop.\ (30.13)]{Knus} for details. Thus, an element of $\HH^1(\Q_p, \CG_\ell)$ can be represented by an equivalence class $[x, y]$ of elements $(x, y) \in \TT(\Q_p)$.

Let  $V$ be a $\ell_\p$-vector space of dimension $n + 1$ and $h$ the hermitian form on $V$ given by the matrix
\[
\left[ 
\begin{array}{cc}
0 & \Id \\
\Id & 0
\end{array}
\right],
\]
where $\Id$ is the identity matrix of dimension $\frac{n+1}{2}$. If $\overline{\G}_\ell$ denotes the adjoint group of $\G_\ell$, then $\overline{\G}_\ell(\Q_p)$ is $\mbox{PGU}(V, h) = \mbox{GU}(V, h) / \ell^\times_p$, where $\mbox{GU}(V, h)$ is the group of unitary similitudes of $(V, h)$. By~\cite[$\S$31.A]{Knus} the connecting homomorphism $\overline{\G}_\ell(\Q_p) \to \HH^1(\Q_p, \CG_\ell)$ has the form
\begin{eqnarray}\label{connecting-for-ajG}
g \mod \ell^\times_\p &\mapsto& [\mu(g), \No(\det(g))],
\end{eqnarray}
where $\mu(g)$ is the multiplier of $g \in \mbox{GU}(V, h)$.

Let $x \in \mu_{n + 1}(\Q_p)$. The element of $g \in \overline{\G}_\ell(\Q_p)$ given by the diagonal matrix 
\[
\mbox{diag}(1, \dots, 1, x, \dots, x) \mod \ell^\times_\p
\]
is clearly contained in a compact torus of $\overline{\G}_\ell(\Q_p)$. Hence it acts trivially on $\Delta_p$ (see~\cite[\S 2.5]{Tits2}), as does its image in $\HH^1(\Q_p, \CG_\ell)$. From~\eqref{connecting-for-ajG} we see that this image is $[x, 1] \in \HH^1(\Q_p, \CG_\ell)$. This element is also obtained as the image of $x$ under the composition $\mu_{n + 1}(\Q_p) \hookrightarrow \TT(\Q_p) \to \HH^1(\Q_p,\CG_\ell)$. We conclude that $\mu_{n + 1}(\Q_p)$ and its quotient $N_p$ act trivially on $\Delta_p$.
\end{pf}



\subsection{Determining $A_\xi$}


By Lemma~\ref{lem:N_p-trivial}, the homomorphism $\xi_p$ defines an action of $W_p$ on $\Delta_p$. We now describe the kernel of this action.

Let $\Ram_\ell$ be the set of finite primes of $\Q$ which ramify in $\ell / \Q$. This is a finite non-empty set. Recall from the introduction that $\# \Ram_\ell = r_\ell$. A nonarchimedean valuation $v$ of $\ell^\times$ is \emph{normalized} if $v(\ell^\times) = \Z$. In the following, for a prime $\mathfrak{p}$ of $\ell$ above $p$, the valuation $v_\mathfrak{p}$ will denote the normalized valuation associated with $\mathfrak{p}$. We denote the $\p$-adic integers of $\ell_\p$ by $\Op$.

Since our group $\G_\ell$ is quasi-split over every $p$-adic field $\Q_p$, we have the following consequence of~\cite[Lemma 2.3 and Prop.\ 2.7]{Borel--Prasad}.


\begin{lem}\label{lem:Borel-Prasad}
Let $p$ be a finite prime with $p \not \in \Ram_\ell$. Then the kernel of the homomorphism $W_p \to \Aut(\Delta_p)$ is given by
\begin{eqnarray}
\Wo_p &=& \ker \left( \prod_{\p|p} \mathcal{O}_\p^\times (\ell_\p^\times)^{n + 1} / (\ell_\p^\times)^{n + 1} \stackrel{\No}{\to} \Q_p^\times / (\Q_p^\times)^{n + 1} \right).
\label{Wo}
\end{eqnarray}
In particular,	$x \in \LL$ acts trivially on $\Delta_p$ if and only if $v_{\mathfrak{p}}(x) \in (n + 1)\Z$ for each prime $\mathfrak{p}$ of $\ell$ above $p$. 
\end{lem}



\begin{rmk}\label{rmk:Wo_p-for-ramified}
For $p \in \Ram_\ell$ we also define $\Wo_p$ as in \eqref{Wo}. As noticed in~\cite[Lemma 5.4]{Borel--Prasad}, for any $p \in \Ram_\ell$ and $x \in \LL$ we have $v_\p(x) \in (n + 1)\Z$. The same argument shows that $W_p = \Wo_p$ for $p \in \Ram_\ell$.
\end{rmk}


We denote by $\LL_{n + 1}$ the subgroup of $\LL$ consisting of elements $x$ with $v(x) \in (n + 1)\Z$ for every normalized  nonarchimedean valuation $v$ of $\ell^\times$. From Lemma~\ref{lem:Borel-Prasad} and Remark~\ref{rmk:Wo_p-for-ramified}, the
sequence~\eqref{A_xi:exact-seq} yields the exact sequence
\begin{eqnarray}\label{A_xi:exact-seq-II}
1 \quad \to \quad A_\xi \quad \to \quad \LL_{n + 1} / (\ell^\times)^{n + 1} \quad \to \quad \prod_{p \in \Ram_\ell} \Aut(\Delta_p).
\end{eqnarray}


\begin{prop}\label{prop:size-of-A_xi}\ 
\begin{enumerate}

\item If $n$ is even, then $A_\xi$ is isomorphic to $\LL_{n+1}/(\ell^\times)^{n+1}$.

\item If $n$ is odd, denote by $\ind$ the index of $A_\xi$ in $\LL_{n + 1} / (\ell^\times)^{n + 1}$. Then $\ind$ is a divisor of $2^{r_\ell}$ with $\ind \ge 2$. In particular, if $r_\ell = 1$ we have $\ind = 2$.

\end{enumerate}
\end{prop}



\begin{pf}
If $n = 2 m$ is even, for $p \in \Ram_\ell$ we have that $\G_\ell / \Q_p$ has $\mbox{C--BC}_m$ as local index, and so $\Aut(\Delta_p) = 1$. The first statement follows immediately from~\eqref{A_xi:exact-seq-II}.

For $n = 2 m - 1$ and $p \in \Ram_\ell$, the local index of $\G_\ell / \Q_p$ is $\mbox{B--C}_m$, for which $\Aut(\Delta_p)$ has order $2$. Then $\ind$ must be a divisor of $2^{r_\ell}$. For $p \in \Ram_\ell$ let $\pi_\p \in \ell^\times$ be any uniformizer with respect to the place $\p$  above $p$. We know from~\cite[$\S$4.2]{Mohammadi--Salehi-Golsefidy} that $\pi_\p \overline{\pi}^{-1}_\p \in \LL$ is a generator of $\Delta_p$ (where $\overline{\pi}_\p$ is the Galois conjugate of $\pi_\p$ in $\ell$). When $p$ is odd, the element $\sqrt{-d}$ is a uniformizer, where $\ell = \Q(\sqrt{-d})$, and we obtain that $-1 \in \LL_{n+1}$ is a generator of $\Aut(\Delta_p)$. This proves that $\ind \ge 2$ when $\Ram_\ell$ contains an odd prime $p$. This is always the case, except when $d$ is $1$ or $2$. For $d=2$, using the same uniformizer $\pi_\p = \sqrt{-d}$ we also see that $\ind = 2$. For $d = 1$ and $\pi_{(2)} = 1 + i$ we have $i \in \LL_{n + 1}$ as generator of $\Delta_2$, and thus $\ind = 2$ as well.
\end{pf}



\subsection{The covolume of $\SU(n, 1) / \Gamma_\ell$}

\label{ss:covolume-of-Gamma_ell}


We denote by $\cI_\ell$ (resp., $\cP_\ell$) the group of fractional (resp., principal) ideals of $\ell$ and by $\cC_\ell = \cI_\ell / \cP_\ell$ the class group of $\ell$. Since the group of integral units in $\LL$ is $\mu(\ell)$, we have the exact sequence (see~\cite[\S 12.6]{Emery} for details)
\begin{eqnarray}\label{ex-seq-LL_mod_ell}
1 \to \mu(\ell) / \mu(\ell)^{n + 1} \to \LL_{n + 1} / (\ell^\times)^{n + 1} \to (\cP_\LL \cap \cI_\ell^{n + 1}) / \cP_\ell^{n + 1} \to 1,
\end{eqnarray}
where $\cP_\LL$ is the subgroup of principal ideals which can be written as $(x)$ with $x \in \LL$.

For an ideal $(x) = \mathfrak{a}^{n + 1} \in \cP_\LL \cap \cI_\ell^{n + 1}$, its ideal norm $\No(x) = [\mathcal{O}_\ell : (x)]$ is obviously in $(\Q^\times)^{n + 1}$. Since $\ell$ is imaginary quadratic this norm is just $\No(x)$. This proves that $\cP_\LL \cap \cI_\ell^{n + 1}$ equals $\cP_\ell \cap \cI_\ell^{n + 1}$.

The quotient of this latter group by $\cP_\ell^{n + 1}$ is isomorphic to the subgroup of $\cC_\ell$ consisting of those elements with order dividing $n + 1$ (see~\cite[proof of Prop.\ 0.12]{Borel--Prasad}). Let $h_{\ell, n + 1}$ be the order of this subgroup. Returning to the sequence~\eqref{ex-seq-LL_mod_ell}, we see that the order of $\LL_{n + 1} / (\ell^\times)^{n + 1}$ is $h_{\ell, n + 1} \cdot \left| \mu(\ell) / \mu(\ell)^{n + 1}\right|$.

The exact sequence
\[
1 \to \mu_{n + 1}(\ell) \to \mu(\ell) \stackrel{x^{n + 1}}{\to} \mu(\ell)^{n + 1} \to 1
\]
implies that $\left| \mu(\ell) / \mu(\ell)^{n + 1}\right|$ equals the order of $\mu_{n + 1}(\ell)$. We conclude that
\begin{eqnarray}\label{bound-LL_n1}
\left| \LL_{n + 1} / (\ell^\times)^{n + 1}\right| & = & h_{\ell, n + 1} \;
\cdot \;
\left| \mu_{n + 1}(\ell) \right|.
\end{eqnarray} 

This equality is the last step in the computation of $[\Gamma_\ell : \Lambda_\ell]$. It follows from the equalities~\eqref{index-Rholfs}, \eqref{index-centre}, \eqref{bound-LL_n1} and Proposition~\ref{prop:size-of-A_xi}
that $[\Gamma_\ell : \Lambda_\ell]$ is equal to $(n + 1) h_{\ell, n + 1}$ if $n$ is even, and $\ind^{-1} (n + 1) h_{\ell, n + 1}$ if $n$ is odd. Together with~\eqref{eq:covol-Lambda}, we get the following.


\begin{prop}\label{prop:minimal-volume-for-each-ell}\ 
\begin{enumerate}

\item If $n$ is even, then $\mu(\SU(n, 1) / \Gamma_\ell)$ equals
\begin{eqnarray*}
\frac{1}{(n + 1) h_{\ell, n + 1}}\;
\Disc_\ell^{\frac{n (n + 3)}{4}} \prod_{j = 1}^n \frac{j!}{(2 \pi)^{j + 1}} \;
\zeta(2) L_\ell(3) \cdots \zeta(n) L_\ell(n+1).
\end{eqnarray*}

\item If $n$ is odd, then $\mu(\SU(n, 1) / \Gamma_\ell)$ is
\begin{eqnarray*}
\frac{\ind}{(n + 1) h_{\ell, n + 1}}\;
\Disc_\ell^{\frac{(n - 1) (n + 2)}{4}} \prod_{j = 1}^n \frac{j!}{(2 \pi)^{j + 1}} \;
\zeta(2) L_\ell(3) \cdots L_\ell(n) \zeta(n + 1),
\end{eqnarray*}
where $\ind$ is as in Proposition~\ref{prop:size-of-A_xi}.
\end{enumerate}
\end{prop}



\subsection{Bounding $[\Gamma_{\ell, \T} : \Lambda_{\ell, \T}]$ and computing $\mvol_\ell(n)$}\label{sec:G-prime}


Suppose that $n$ is odd, and consider the lattice $\Gamma_{\ell,\T}$ (see $\S$\ref{sec:constuction-of-Lambda}). The index $[\Gamma_{\ell, \T} : \Lambda_{\ell, \T}]$ can be expressed by the same formula as in~\eqref{index-Rholfs}, where $A \subset \HH^1(\Q, \CG_\ell)$ now acts on the local Dynkin diagrams $\Delta_p$ associated with $\G_{\ell, \T}$.

Lemma~\ref{lem:Borel-Prasad} is not valid for primes $p$ where $\G_{\ell, \T} / \Q_p$ fails to be quasi-split. In particular $A_\xi$ does not need to be a subgroup of $\LL_{n + 1} / (\ell^\times)^{n + 1}$. However, $\G_{\ell, \T} / \Q_p$ is not quasi-split only for a finite number of primes $p$, so $A_\xi$ will be commensurable with $\LL_{n + 1} / (\ell^\times)^{n + 1}$, and it is still possible to bound $[\Gamma_{\ell, \T} : \Lambda_{\ell, \T}]$. A careful analysis (see~\cite[$\S$12.5]{Emery} for details) shows that here we have the following bound.


\begin{lem}\label{lem:bound-for-index}
For any $\Lambda_{\ell, \T}$,
\begin{eqnarray}\label{bound-for-index}
[\Gamma_{\ell, \T} : \Lambda_{\ell, \T}] &\le& (n + 1)^{\#\hat{\T} + 1} h_{\ell, n + 1},
\end{eqnarray}
where $\hat{\T}$ denotes the set of unramified primes $p$ where $\G_{\ell, \T} / \Q_p$ is not quasi-split.
\end{lem}


Clearly $\hat{\T} \subset \T$, where $\T$ is the set of primes defined
in $\S$\ref{sec:constuction-of-Lambda}. Note that the volume of
$\SU(n, 1) / \Gamma_{\ell, \T}$ then has a factor of $\lambda_p / (n +
1)$ for every $p \in \hat \T$, where the factors $\lambda_p$ are given
in Equation \eqref{eq:lambda-factor}. One can compute the factors $\lambda_p$ from~\cite[Table 1]{Ono} and see that $\lambda_p / (n + 1) > 1$ for every $p$ and every $n \geq 3$. It follows then from Lemma~\ref{lem:bound-for-index} and the explicit expression for $[\Gamma_\ell : \Lambda_\ell]$ that the covolume of $\Gamma_{\ell}$ is always smaller than the covolume of $\Gamma_{\ell,\T}$.

We are now ready for the following.


\begin{pf}[Computation of $\mvol_\ell(n)$ for Theorems~\ref{thm:even} and~\ref{thm:odd}]
Let $\Gamma < \PU(n, 1)$ be a lattice of minimal covolume for its commensurability class. It follows from the above that $\Gamma$ is the image in $\PU(n, 1)$ of $\Gamma_\ell$.

As in~\cite[$\S$4]{Borel--Prasad}, we have that the Euler--Poincar\'e
measure on $\SU(n,1)$ equals 
\[
 \chi(\C \mathbb P^n) \mu = (n + 1) \mu \; ,
\]
since $\C \mathbb P^n$ is the compact dual to $\Hysp_\C^n$. Since $\SU(n, 1)$ covers $\PU(n, 1)$ with degree $n + 1$, we see that
\begin{eqnarray*}
	\EP(\PU(n,1)/\Gamma) &=&  (n + 1)^2 \mu(\SU(n, 1) / \Gamma_\ell).
\end{eqnarray*}
The equations for $\mvol_\ell(n)$ in Theorems \ref{thm:even} and \ref{thm:odd} now follow from Proposition \ref{prop:minimal-volume-for-each-ell} followed by the standard functional equations for $\zeta$ and $L_\ell$ (see \cite[$\S$15]{Swinnerton-Dyer}).
\end{pf}



\subsection{The proof of Theorem~\ref{thm:smallest volume}}


In this section, we prove the first part of Theorem~\ref{thm:smallest volume}. The case $n = 2$ appears in~\cite{Stover}.


\begin{prop}\label{nonuniform minimal}
For all $n \geq 2$ the minimal covolume nonuniform arithmetic lattice in $\SU(n, 1)$ is $\Gamma_\ell$ for $\ell = \Q(\sqrt{-3})$ and $\T = \emptyset$.
\end{prop}



\begin{pf}
Suppose that the lattice of smallest covolume comes from some other imaginary quadratic field $\ell$. By $\S$\ref{sec:G-prime}, this lattice is some $\Gamma_\ell$. Let $\ell_0 = \Q(\sqrt{-3})$.


We now give an upper bound for the discriminant of $\ell$. First, set $s = n(n + 3) / 4$ for $n$ even and $s = (n - 1)(n + 2) / 4$ for $n$ odd. Then, assume that
\[
\frac{1}{[\Lambda_\ell : \Gamma_\ell]} \Disc_\ell^s \prod_{j = 1}^n \frac{j!}{(2 \pi)^{j + 1}} \;
\zeta(2) L_\ell(3) \cdots \left\{ \begin{array}{lr}
\hspace{-0.2cm} \zeta(n + 1) & \mbox{$n$ odd}\\
\hspace{-0.2cm} L_\ell(n + 1) & \mbox{$n$ even}
\end{array}\right.
\]
\[
\leq \frac{q^\prime}{n + 1} 3^s \prod_{j = 1}^n \frac{j!}{(2 \pi)^{j + 1}} \;
\zeta(2) L_{\ell_0}(3) \cdots \left\{ \begin{array}{lr}
\hspace{-0.2cm} \zeta(n + 1) & \mbox{$n$ odd}\\
\hspace{-0.2cm} L_{\ell_0}(n + 1) & \mbox{$n$ even}
\end{array} \right.
\]
When $n$ is odd, we can cancel the factor of $\zeta(n + 1)$ from each side. For all $n$, we cancel the factors of $j! / (2 \pi)^{j + 1}$ to get
\begin{equation}\label{first nonuniform estimate}
\frac{\Disc_\ell^s}{[\Lambda_\ell : \Gamma_\ell]} \zeta(2) L_\ell(3) \cdots L_\ell(n + t) \leq \frac{3^s q^\prime}{n + 1} \zeta(2) L_{\ell_0}(3) \cdots L_{\ell_0}(n + t),
\end{equation}
where $t = 0$ for $n$ odd and $t = 1$ for $n$ even.

Comparing terms of each series expansion, we see that $\zeta(k) L_{\ell_0}(k + 1) \leq \zeta(2) L_{\ell_0}(3)$ for any $k \geq 2$. Therefore, we have
\[
\frac{3^s q^\prime}{n + 1} \zeta(2) L_{\ell_0}(3) \cdots L_{\ell_0}(n + t) \leq \frac{3^s q^\prime}{n + 1} (\zeta(2) L_{\ell_0}(3))^{\frac{n - t}{2}}
\]
\[
= \frac{3^s q^\prime}{n + 1} \left( \frac{2 \pi ^5}{3^{11 / 2}} \right)^{\frac{n - t}{2}} \leq 2 \frac{3^s}{n + 1} \left( \frac{2 \pi ^5}{3^{11 / 2}} \right)^{\frac{n - t}{2}}.
\]
One the other hand, from Proposition~\ref{prop:minimal-volume-for-each-ell}, we have
\[
\frac{\Disc_\ell^s}{[\Lambda_\ell : \Gamma_\ell]} \zeta(2) L_\ell(3) \cdots L_\ell(n + t) \geq
\]
\[
\frac{\Disc_\ell^s}{(n + 1) h_{\ell, n + 1}} \zeta(2) L_\ell(3) \cdots L_\ell(n + t).
\]

Then we have $h_{\ell, n + 1} \leq h_\ell$, and a key ingredient of the proof of the Brauer--Siegel Theorem implies that
\[
h_\ell \leq \mu_\ell m (m - 1) (m - 1)! \left( \frac{D_\ell}{4 \pi^2} \right)^{m / 2} \zeta_\ell(m)
\]
for any integer $m > 1$ and imaginary quadratic field $\ell$, where $\mu_\ell$ is the order of the group of roots of unity in $\ell$ and $\zeta_\ell$ is the zeta function of $\ell$. See the proof of Lemma 1 in $\S$1 of Chapter XVI of \cite{Lang}. For $\ell$ not $\Q(i)$ or $\Q(\sqrt{-3})$, $\mu_\ell = 2$.

From here forward assume $D_\ell > 4$, so
\[
h_\ell \leq 2 m (m - 1) (m - 1)! \left( \frac{D_\ell}{4 \pi^2} \right)^{m / 2} \zeta_\ell(m)
\]
for all $m > 1$. Taking $m = 2$, we get
\[
\frac{\Disc_\ell^s}{[\Lambda_\ell : \Gamma_\ell]} \zeta(2) L_\ell(3) \cdots L_\ell(n + t)
\]
\[
\geq \frac{D_\ell^{s - 1} \pi^2}{\zeta_\ell(2)} \zeta(2) L_\ell(3) \cdots L_\ell(n + t).
\]
Then $\zeta(2 k) L_\ell(2 k + 1) > 1$ for all $k > 1$ and $\zeta_\ell(2) \leq \zeta(2) = \pi^2 / 6$, so we have
\[
\frac{\Disc_\ell^s}{[\Lambda_\ell : \Gamma_\ell]} \zeta(2) L_\ell(3) \cdots L_\ell(n + t) \geq 6 D_\ell^{s - 1}.
\]

We conclude that if $\ell$ gives smaller covolume than $\Q(\sqrt{-3})$, then
\[
D_\ell \leq 3 \left(\frac{2 \pi ^5}{3^{11 / 2}} \right)^{\frac{n - t}{2(s - 1)}}.
\]
An immediate computation shows that this is bounded above by $4$ for all $n \geq 4$, so we need only consider $\Q(i)$. For $\Q(i)$, we replace $\mu_\ell$ with 4 in the above computations and see that $\Q(i)$ must give larger covolume than $\Q(\sqrt{-3})$ as well, provided that $n \geq 5$. One then checks by hand that the smallest covolume comes from $\Q(\sqrt{-3})$ for $n = 2, 3, 4$. This completes the proof of the proposition.
\end{pf}



%% file: multiplicities6.tex

In this section we address the question of determining the number of isomorphism classes of arithmetic subgroups of $\G_\ell$ of minimal covolume.  In \S\ref{ss:expressions-for-n_d}--\ref{ss:class-number-by-ideles} we will first count conjugacy classes of arithmetic subgroups of $\G_\ell(\R)$, and then in \S \ref{ss:count-up-isom} we show how this implies our results on counting isomorphism classes stated in $\S$\ref{intro}.  
\subsection{Counting conjugacy classes}\label{ss:expressions-for-n_d}


Let us fix an identification $\G_\ell(\R) = \SU(n,1)$.  The arithmetic subgroup $\Lambda_\ell$ = $\Lambda_\ell(P)$ and its normalizer $\Gamma_\ell = \Gamma_\ell(P)$ in $\SU(n, 1)$ depend on the choice of a coherent collection $P \in \ccol_\ell$ (see $\S$\ref{sec:constuction-of-Lambda}). Though the covolume does not depend on this choice, two different coherent collections of parahorics with the same local types may determine lattices that are not conjugate under the action of $\PU(n,1)$ by inner automorphisms. 

In this section it will be useful to consider $\Lambda_\ell$ and $\Gamma_\ell$ as functions 
\begin{eqnarray*}
\Lambda_\ell &:& \ccol_\ell \to \left\{ \mbox{lattices in} \G_\ell(\Q) \right\}, \\
\Gamma_\ell &:& \ccol_\ell \to \left\{ \mbox{lattices in} \G_\ell(\R) \right\}.
\end{eqnarray*}
Then $\PU(n, 1)$ acts by conjugation on both $\Lambda_\ell(\ccol_\ell)$ and $\Gamma_\ell(\ccol_\ell)$.
Let $\bG_\ell$ be the adjoint group of $\G_\ell$. For any field extension $K / \Q$ the group $\bG_\ell(K)$ corresponds to the group of inner automorphisms of $\G_\ell$ defined over $K$. For $K = \R$ we have $\bG_\ell(\R) = \PU(n, 1)$. The group $\bG_\ell(\Q)$ acts on the set $\ccol_\ell$ diagonally by conjugating each component $P_p \subset \G_\ell(\Q_p)$ of a coherent collection $P = (P_p)_{p \in \Vf} \in \ccol_\ell$.  

\begin{lem}\label{Gamma-conj-cl=Lambda-conj-cl}
The number $\nbo_\ell$  of $\PU(n, 1)$-conjugacy classes in $\Gamma_\ell(\ccol_\ell)$ is equal to the number of $\bG_\ell(\Q)$-conjugacy classes in $\ccol_\ell$.
\end{lem}



\begin{pf}
Since arithmetic subgroups are Zariski dense in the $\Q$-group $\G_\ell$ they can only be conjugated by inner automorphisms that are in $\bG_\ell(\Q)$. Thus we must count the number of $\bG_\ell(\Q)$-conjugacy classes in $\Gamma_\ell(\ccol_\ell)$. Let $P$ and $P'$ be two coherent collection in $\ccol_\ell$. Since $\Gamma_\ell$ is by definition the normalizer of $\Lambda_\ell$, if $\Lambda_\ell(P)$ and $\Lambda_\ell(P')$ are conjugate, then $\Gamma_\ell(P)$ and $\Gamma_\ell(P')$ are also conjugate by the same element. The converse is also true and follows from maximality of $\Lambda_\ell(P)$ and $\Lambda_\ell(P')$ in $\G_\ell(\Q)$.

This proves that $\nbo_\ell$ is equal to the number of $\bG_\ell(\Q)$-conjugacy classes in $\Lambda_\ell(\ccol_\ell)$. By Strong Approximation, every coherent collection $P$ can be recovered by the principal arithmetic subgroup $\Lambda_\ell(P)$. In particular, if $\Lambda_\ell(P)$ and $\Lambda_\ell(P')$ are conjugate, then $P$ and $P'$ are conjugate. This finishes the proof.
\end{pf}


Choose $P = (P_p) \in \ccol_\ell$. The adelic group $\bG_\ell(\Af)$  acts on $\ccol_\ell$, each component $g_p$ of $g = (g_p)_{p \in \Vf} \in \bG_\ell(\Af)$ acting by conjugation on $P_p$. The group $\bG_\ell(\Q)$ is diagonally embedded in $\bG_\ell(\Af)$ as usual. Then the action of the adelic group restricted to $\bG_\ell(\Q)$ coincides with the action of $\bG_\ell(\Q)$ on $\ccol_\ell$ as above.

Fix a coherent collection $P \in \ccol_\ell$.  For each $p \in \Vf$ we denote by $\bP_p$ the stabilizer of $P_p$ in $\bG_\ell(\Q_p)$. We can identify the $\bG_\ell(\Af)$-conjugacy class of $P$ with the set 
$\bG_\ell(\Af) / \prod_{p} \bP_p $. The $\bG_\ell(\Af)$-conjugacy class of $P$ in $\ccol_\ell$ is divided into $\bG_\ell(\Q)$-conjugacy classes, which are in bijection with the double cosets
\begin{eqnarray}\label{class-set-of-group}
\bG_\ell(\Q) \backslash \bG_\ell(\Af) / \prod_{p \in \Vf} \bP_p.
\end{eqnarray}
The cardinality of this set is called the \emph{class number of $\bG_\ell$ relative to $P$}, and is known to be finite (see for example~\cite[Prop.\ 3.9]{Borel--Prasad}).


\subsection{Computing the class number}\label{ss:class-number-indep-of-P}


For each prime $p$, consider the connecting homomorphism \begin{eqnarray}\label{connecting-delta_p} \delta_p &:& \bG_\ell(\Q_p) \to \HH^1(\Q_p, \CG_\ell), \end{eqnarray} where, as in $\S$\ref{nonuniform prasad}, the group $\CG_\ell$ denotes the center of $\G_\ell$. Let $\prod_{p \in \Vf}' \HH^1(\Q_p, \CG_\ell)$ be the restricted direct product with respect to the groups $\delta_p(\bP_p)$. The group $\HH^1(\Q, \CG_\ell)$ is diagonally embedded in $\prod'_{p \in \Vf} \HH^1(\Q_p, \CG_\ell)$. Following~\cite[\S 5.2]{Prasad--Yeung1} or~\cite[\S 6.1]{Belolipetsky--Emery} we know that there is a bijection between the double coset~\eqref{class-set-of-group} and the group
\begin{eqnarray}\label{class-group-by-cohomology}
\CC_\ell(P) &=& \frac{\prod'_{p \in \Vf} \HH^1(\Q_p,\CG_\ell)}{\delta\left( \bG_\ell(\Q) \right) \cdot \prod_p \delta_p(\bP_p)},
\end{eqnarray}
where $\delta$ is the connecting homomorphism from $\bG_\ell(\Q)$ to $\HH^1(\Q, \CG_\ell)$. Note that the quotient on the right is an abelian group.

Since $P$ is a coherent collection in $\ccol_\ell$ we have that each $P_p$ has special type. An automorphism of the local Dynkin diagram $\Delta_p$ of $\G_\ell(\Q_p)$ that fixes a special type  is necessarily trivial~\cite[Prop.\ 12.2]{Emery}.  It follows that for each prime $p$ the group $\delta_p(\bP_p)$ is precisely the kernel $\HH^1(\Q_p, \CG_\ell)_{\xi_p}$ of the map $\xi_p$ introduced in $\S$\ref{index-comput}. In particular, this kernel is independent of $P_p$. Thus the order of $\CC_\ell = \CC_\ell(P)$ does not depend on $P \in \ccol_\ell$. We will denote this order by $\cln_\ell$ and call it the \emph{class number} of $\bG_\ell$.  

\begin{prop}\label{prop:n_p-in-fct-of-h_d}
Let $r_\ell$ be the number of ramified primes in $\ell / \Q$. Then
\begin{eqnarray*}
\nbo_\ell &=& \left\{ \begin{array}[h]{cl}
2^{r_\ell} \cln_\ell & \mbox{ if $n$ is even};\\
\cln_\ell & \mbox{ if $n$ is odd}.
\end{array} \right.
\end{eqnarray*}
\end{prop}



\begin{pf}
By Lemma~\ref{Gamma-conj-cl=Lambda-conj-cl}, $\nbo_\ell$ is given by the number of $\bG_\ell(\Q)$-conjugacy classes in the set $\ccol_\ell$. By the above discussion, the set of these classes is divided into $\bG_\ell(\Af)$-conjugacy classes, where each $\bG_\ell(\Af)$-conjugacy classes identifies exactly $\cln_\ell$ of the $\bG_\ell(\Q)$-conjugacy classes. Therefore, we only need to determine the number of $\bG_\ell(\Af)$-conjugacy classes in $\ccol_\ell$. If $p \not \in \Ram_\ell$, the type of $P_p$ for $P \in \ccol_\ell$ is hyperspecial and $\bG_\ell(\Q_p)$ acts transitively on hyperspecial types~\cite[\S 2.5]{Tits2}. If $n$ is odd and $p \in \Ram_\ell$ then $\bG_\ell(\Q_p)$ permutes the two special points of the local Dynkin diagram $\Delta_p$. It follows that all coherent collections of $\ccol_\ell$ are conjugate in this case.

When $n$ is even we have $\Aut(\Delta_p) = 1$ for $p \in \Ram_\ell$ and, in particular, the two special points of $\Delta_p$ are not conjugate. This shows that for even $n$ there are as many $\bG_\ell(\Af)$-conjugacy classes in $\ccol_\ell$ as there are possible choices of a special point in each $\Delta_p$ for $p \in \Ram_\ell$. That is, there are $2^{r_\ell}$ $\bG_\ell(\Af)$-conjugacy classes.
\end{pf}




\subsection{More on the class number}\label{ss:class-number-as-W_p}


Following $\S$\ref{ss:class-number-indep-of-P} we know that $\cln_\ell$ is given by the size of the group
\begin{eqnarray*}\label{h_d-as-size-HH1}
\CC_\ell &=&  \frac{\prod'_{p \in \Vf} \HH^1(\Q_p, \CG_\ell)}{ \delta\left( \bG_\ell(\Q) \right) \cdot \prod_p \HH^1(\Q_p, \CG_\ell)_{\xi_p}}.
\end{eqnarray*}
Using Lemma~\ref{lem:N_p-trivial}, we can factor out the product $\prod_p N_p$ from the numerator and denominator. With the description of $\delta\left( \bG_\ell(\Q) \right)$ in $\S$\ref{ss:class-number-indep-of-P} we obtain an isomorphism 
\begin{eqnarray}\label{h_d-as-size-W}
\CC_\ell &\cong& \frac{\prod'_{p \in \Vf} W_p}{\LL/(\ell^\times)^{n + 1} \cdot \prod_p W_{\xi_p}},
\end{eqnarray}
where $W_{\xi_p}$ is the kernel of the homomorphism $W_p \to \Aut(\Delta_p)$ induced by $\xi_p$. 

For $p \not \in \Ram_\ell$ we have $W_{\xi_p} = \Wo_p$ (cf., Lemma~\ref{lem:Borel-Prasad}). If the dimension $n$ is even we have $\Aut(\Delta_p) = 1$ for $p \in \Ram_\ell$. Using Remark~\ref{rmk:Wo_p-for-ramified} we see that $W_{\xi_p} = \Wo_p$ in this case as well, so
\begin{eqnarray}\label{h_d-as-size-W-Wo}
\CC_\ell &\cong& \frac{\prod'_{p \in \Vf} W_p}{ \LL/(\ell^\times)^{n + 1} \cdot \prod_p \Wo_p}.
\end{eqnarray}

Now suppose that $n$ is odd. For $p \in \Ram_\ell$ the subgroup $W_{\xi_p}$ is of index two in $W_p = \Wo_p$. However, if $\# \Ram_\ell  = 1$ the proof of Proposition~\ref{prop:size-of-A_xi} shows that 
\begin{eqnarray*}
\LL/(\ell^\times)^{n+1} \cdot \prod_p W_{\xi_p} &=& \LL/(\ell^\times)^{n+1} \cdot \prod_p \Wo_p.
\end{eqnarray*}
Therefore, assuming $r_\ell = 1$, Equation \eqref{h_d-as-size-W-Wo}
still holds.

\subsection{Idelic formulation}\label{ss:class-number-by-ideles}


Let $\Idl_\ell$ be the group of finite ideles of $\ell$.  We denote by $\Idlo_\ell$ the subgroup of integral finite ideles~\cite[$\S$9]{Swinnerton-Dyer}. The norm map $\No$ extends to ideles as a map  $\No: \Idl_\ell \to \Idl_\Q$, defined component-wise, where $\Idl_\Q$ is the group of ideles of $\Q$. Consider the kernel 
\begin{eqnarray}\label{Idl_LL}
\Idl_\LL &=&  \ker\left( \Idl_\ell \stackrel{\No}{\to} \Idl_\Q / \Idl_\Q^{n + 1} \right).
\end{eqnarray}

Equation \eqref{h_d-as-size-W-Wo}, which holds when $n$ is even or $r_\ell = 1$, can be rewritten using ideles as
\begin{eqnarray}
\CC_\ell &\cong& \frac{\Idl_\LL / \Idl_\ell^{n + 1}}{ \LL/(\ell^\times)^{n + 1} \cdot \Idlo_\ell \Idl_\ell^{n + 1} / \Idl_\ell^{n + 1}} \nonumber \\
& \cong &  \Idl_\LL / \left( \LL \Idlo_\LL \Idl_\ell^{n + 1} \right),
\end{eqnarray}
with $\LL$ embedded diagonally in $\Idl_\LL$, and $\Idlo_\LL = \Idlo_\ell \cap \Idl_\LL$.

Recall $h_{\ell, n + 1}$ is the order of the subgroup of the class group $\cC_\ell$ consisting of those elements with order dividing $n + 1$. We then have the following.


\begin{lem}\label{lem:order-of-I_LL-mod}
The order of the group $\Idl_\LL / \left( \LL \Idlo_\LL \Idl_\ell^{n + 1} \right)$ divides $2 h_{\ell, n + 1}$. If $8 \nmid n + 1$, this order divides $h_{\ell, n + 1}$. 
\end{lem}



\begin{pf}
By definition of $\LL$ and $\Idl_\LL$, we have the exact sequence
\begin{eqnarray*}
1 \quad \to \quad \LL \quad \to \quad \ell^\times \cap \Idl_\LL \quad \stackrel{\No}{\to} \quad \Q^\times / (\Q^\times)^{n + 1} \quad \stackrel{\phi}{\to} \quad \prod_{p \in \Vf} \Q_p^\times / (\Q_p^\times)^{n + 1},
\end{eqnarray*}
where $\phi$ is the natural diagonal map. This shows that the index of $\LL$ in $\ell^\times \cap \Idl_\LL$ has the same size as the kernel of $\phi$. According to~\cite[Thm 9.1.11]{NSW} this kernel has order at most $2$, and is trivial if $8 \nmid n + 1$.

Each idele of $\ell$ determines a fractional ideal of $\ell$ by a standard procedure described, for example, in~\cite[$\S$9]{Swinnerton-Dyer}. This induces a homomorphism
\begin{eqnarray}\label{ideles-to-ideals}
\Idl_\LL &\to& \cI_\ell / \cP_\ell  \cI_\ell^{n + 1},
\end{eqnarray}
whose kernel is $(\ell^\times \cap \Idl_\LL) \Idlo_\LL \Idl_\ell^{n + 1}$. Now $\cI_\ell / \cP_\ell \cI_\ell^{n + 1}$ is the quotient $\cC_\ell / \cC_\ell^{n + 1}$ of the class group $\cC_\ell$, whose order is easily seen to be $h_{\ell, n + 1}$. Together with the index of $\LL$ in $\ell^\times \cap \Idl_\LL$ this concludes the proof.
\end{pf}


We sum up the results about $\cln_\ell$ in the following proposition. 


\begin{prop}\label{prop:value-of-h_d}\ 
\begin{enumerate}

\item If $n$ is even, then $\cln_\ell$ is a divisor of $h_{\ell, n + 1}$

\item Suppose that $n$ is odd and $r_\ell = 1$. Then $\cln_\ell$ is a divisor of $2 h_{\ell, n + 1}$. Moreover, if $8 \nmid n + 1$, then $\cln_\ell$ is a divisor of $h_{\ell, n + 1}$. 

\end{enumerate}
\end{prop}


\subsection{Counting up to isomorphism}\label{ss:count-up-isom}

Any nonuniform lattice in $\PU(n,1)$ associated with $\ell/\Q$ and of minimal covolume among them corresponds to the projection of a subgroup $\Gamma_\ell \subset \G_\ell(\R)$.  It is clear that the number $\nb_\ell = \nb_\ell(n)$ of isomorphism classes of such lattices is at most equal to the number $\nbo_\ell$ of $\PU(n,1)$-conjugacy classes. Since $\PU(n,1)$ is of index $2$ in its automorphism group, Mostow-Prasad rigidity shows that $\nb_\ell$ is at least equal to $\nbo_\ell/2$. However, when $n$ is even and $p \in \Ram_\ell$, the two special vertices on $\Delta_p$ correspond to nonisomorphic parahoric subgroups in $\G_\ell(\Q_p)$, and thus different choices at this prime will give nonisomorphic lattices. This shows that for $n$ even $\nb_\ell$ is at least equal to $2^{r_\ell}$. From Propositions \ref{prop:n_p-in-fct-of-h_d} and \ref{prop:value-of-h_d} we thus obtain
the following result.

\begin{prop}\label{prop:value-of-nb}\ 
\begin{enumerate}

	\item If $n$ is even, then $\nb_\ell \geq
		2^{r_\ell}$ and $2^{r_\ell-1} \cln_\ell \leq \nb_\ell
		\leq 2^{r_\ell} \cln_\ell$, where $\cln_\ell$ is a
		divisor of $h_{\ell,n+1}$.

	\item Suppose that $n$ is odd and $r_\ell = 1$.
		Then $\cln_\ell/2 \leq \nb_\ell \leq \cln_\ell$, where
	$\cln_\ell$ is a divisor of $2 h_{\ell, n + 1}$. 
	Moreover, if $8 \nmid n + 1$, then $\cln_\ell$ is a divisor of $h_{\ell, n + 1}$. 

\end{enumerate}
\end{prop}



%% file: growth6.tex
\subsection{Overall minimal covolume}\label{sec:overall}

We now prove Theorem \ref{thm:smallest overall} assuming Theorem \ref{thm:volume growth}, which we prove in the sequel. By Theorem~\ref{thm:smallest volume}, a nonuniform arithmetic lattice $\Lambda < \PU(m,1)$ such that $|\chi(\Gamma)| > |\chi(\Lambda)|$ for any $\Gamma < \PU(n,1)$ is necessarily associated with the extension $\Q(\sqrt{-3})/\Q$.  Using the formulas for $\nu_{\Q(\sqrt{-3})}(n)$ given in Theorem~\ref{thm:even} and~\ref{thm:odd}, and the fact that for $n$ high enough the covolume grows (Theorem~\ref{thm:volume growth}), we find that the minimum appears for $m = 9$, where the lattice $\Gamma_0$ of minimal covolume is unique.

Equip complex hyperbolic space with the metric under which holomorphic sectional curvatures are $-1$. We apply Chern--Gauss--Bonnet (see \cite{Hersonsky--Paulin, Parker}) to obtain
\begin{eqnarray}
\vol(\mathbf H_\mathbb C^n / \Gamma) &=&  \frac{(-4 \pi)^n}{(n + 1)!} \; \chi(\Gamma).
	\label{eq:Chern-Gauss-Bonnet}
\end{eqnarray}
Repeating the computations for the Euler characteristic with the complex hyperbolic volume, one sees that replacing $\EP$ by $\vol$ does not change the fact that the smallest volume appears in dimension $m = 9$.  This completes the proof of Theorem~\ref{thm:smallest overall}.


\subsection{Asymptotic behavior of the minimal covolume}
\label{ss:asymptotic}

We proceed with the proof of Theorem~\ref{thm:volume growth}.


\begin{pf}[Proof of Theorem~\ref{thm:volume growth}]
Let $\ell/\Q$ be an imaginary quadratic extension. The ratio $\q_\ell(n) = \nu_\ell(n+1)/\nu_\ell(n)$ is given by:
\begin{eqnarray*}
	\q_\ell(n) &=&  \Disc_\ell^{\pm \frac{n+1}{2}} \epsilon_\ell(n)^{\mp 1} \frac{n + 2}{n + 1} \frac{(n + 1)!}{(2 \pi)^{n + 2}} F(n + 2),
\end{eqnarray*}
where $F$ is $L_{\ell}$ when $n$ is odd and is $\zeta$ when $n$ is even, and the signs $\pm$ depend on the parity of $n$. Since $\ell$ is fixed, $\epsilon_\ell(n)$ is bounded. We also have that
\begin{eqnarray*}
\lim_{n \to \infty}  \frac{n + 2}{n + 1} F(n + 2) &=&  1.
\end{eqnarray*}
Therefore, for $n$ large and up to a constant, $\q_\ell(n)$ grows like 
\begin{eqnarray*}
	\left(\frac{\Disc_\ell^{\pm 1/2}}{(2 \pi)^{n + 2}}\right)^{n+2}
	(n+1)!\;,
\end{eqnarray*}
which is super-exponential.
\end{pf}

Note that the analogous statement for complex hyperbolic volumes is immediate, since the ratio of the additional factors given in~\eqref{eq:Chern-Gauss-Bonnet} decreases linearly in $n$. It also follows that the minimal volume of a noncompact arithmetic complex hyperbolic $n$-manifold grows super-exponentially in $n$. We now analyze the difference between our results and other estimates for volumes of noncompact complex hyperbolic orbifolds.

A lower bound for the volume of a complex hyperbolic $n$-manifold with $k$ cusps, due to Hwang~\cite{Hwang}, is
\[
k \frac{(4 \pi)^n}{n! (P(4) - P(2))} \left( 1 - \frac{n + 1}{P(4) - P(2)} \right),
\]
where
\[
P(\ell) = \frac{(n \ell + n + \ell)!}{n! (n \ell + \ell)!}.
\]
This bound comes from algebraic geometry, and requires smoothness. This bound improves upon earlier estimates due to Parker~\cite{Parker} and Hersonsky--Paulin~\cite{Hersonsky--Paulin}. Note that if $k$ is bounded, this estimate decays to $0$ as $n \to \infty$. Comparison with Theorem~\ref{thm:volume growth} raises the question of how the number of cusps grows as $n \to \infty$.

More recently, generalizing work of Adeboye in real hyperbolic space~\cite{Adeboye}, Fu, Li, and Wang~\cite{Fu--Li--Wang} gave a lower bound for the volume of a complex hyperbolic $n$-orbifold depending only on the largest order of a torsion element of its fundamental group. Adeboye's work has also been extended to complex hyperbolic space by work in progress by Adeboye--Wei~\cite{Adeboye--Wei}. These bounds also go to zero as $n$ goes to infinity. It would be interesting to better understand the gap between these geometric volume bounds and our arithmetic results.